\documentstyle[11pt]{article}
\newlength{\standardunitlength}
 \newtheorem{lemma}{Lemma}
\newtheorem{theorem}{Theorem} 
\newenvironment{proof}{\noindent {\sc Proof:}}{$\Box$ \vspace{2 ex}}

\begin{document}

\begin{center}
Orbifold Euler characteristics and the number of commuting
$m$-tuples in the symmetric groups
\end{center}

\begin{center}
By Jim Bryan$^{*}$ and Jason Fulman
\end{center}

\begin{center}
University of California at Berkeley and Dartmouth College
\end{center}

\begin{center}
Direct Correspondence to:
\end{center}

\begin{center}
Jason Fulman
\end{center}

\begin{center}
Dartmouth College
\end{center}

\begin{center}
Department of Mathematics
\end{center}

\begin{center}
6188 Bradley Hall
\end{center}

\begin{center}
Hanover, NH 03755
\end{center}

\begin{center}
email:jason.e.fulman@dartmouth.edu
\end{center}

\begin{center}
\end{center}

\begin{center}
\end{center}

\begin{center}
\end{center}

\begin{center}
\end{center}

\begin{center}
\end{center}

\begin{center}
\end{center}

$^*$ Supported in part by a grant from the Ford Foundation.

\newpage
\begin{abstract}
Generating functions for the number of commuting $m$-tuples in the
symmetric groups are obtained. We define a natural sequence of ``orbifold
Euler characteristics'' for a finite group $G$ acting on a manifold $X$.
Our definition generalizes the ordinary Euler characteristic of $X/G $ and
the string-theoretic orbifold Euler characteristic. Our formulae for
commuting $m$-tuples underlie formulas that generalize the results of
Macdonald and
Hirzebruch-H\"ofer concerning the ordinary and string-theoretic Euler
characteristics of symmetric products.
\end{abstract}

\section{Introduction}
Let $X$ be a manifold with the action of a finite group $G$. The Euler
characteristic of the quotient space $X/G$ can be computed by the Lefschetz
fixed
point formula:
$$
\chi (X/G)=\frac{1}{|G|}\sum_{g\in G}\chi (X^{g})
$$
where $X^{g}$ is the fixed point set of $g$. Motivated by string theory,
physicists have defined an ``orbifold characteristic'' by
$$
\chi (X,G)=\frac{1}{|G|}\sum_{gh=hg}\chi (X^{(g,h)})
$$
where the sum runs over commuting pairs and $X^{(g,h)}$ denotes the common
fixed point set of $g$ and $h$.

We introduce a natural sequence
of orbifold Euler characteristics $\chi _{m}(X,G)$ for $m=1,2,\ldots$ so
that $\chi (X/G)$ and $\chi (X,G)$ appear as the first two terms. Namely,
if we denote by $Com (G,m)$ the set of mutually commuting $m$-tuples
$(g_{1},\ldots,g_{m})$ and by $X^{(g_{1},\ldots,g_{m})}$ the simultaneous
fixed point set, then we define the {\em $m$-th orbifold characteristic} to
be
\begin{equation}\label{eqn: mth orb euler char defn}
\chi _{m}(X,G)=\frac{1}{|G|}\sum_{Com (G,m)}\chi (X^{(g_{1},\ldots,g_{m})}).
\end{equation}
In the case of a symmetric product, {\em i.e.  } $X$ is the $n$-fold
product $M^{n}$ and $G$ is the symmetric group $S_{n}$, there are
combinatorial formulas for $\chi _{1}$ and $\chi _{2}$ due to Macdonald
\cite{Mac} and Hirzebruch-H\"ofer \cite{Hi-Ho} respectively. The main
result of this note (Theorem \ref{mainresult}) is a generalization of those
formulas to $\chi _{m}$
for arbitrary $m$. In the case where $M$ has (ordinary) Euler
characteristic 1, our formulas specialize to generating functions for
$|Com (S_{n},m)|$, the number of commuting $m$-tuples in $S_{n}$.

        Finally, we remark that the first two terms in our sequence $\chi
_{m}(X,G)$ of
orbifold Euler characteristics are the Euler characteristics of the
cohomology theories $H^{*}_{G}(X;\mathbf{Q})$ and
$K^{*}_{G}(X;\mathbf{Q})$ respectively. This was observed by Segal,
\cite{Atiyah-Segal} who was led to speculate that the
heirarchy of generalized cohomology theories investigated by Hopkins and
Kuhn \cite{Hopkins-Kuhn} may have something to do with the sequence
of Euler characteristics defined in this paper (our definition is
implicitly suggested in \cite{Atiyah-Segal}).
We hope that our combinatorial formulas will provide clues to the nature
of these theories.

\section{Formulae}

In this section we specialize to the case of symmetric products so that
$X=M^{n}$ and $G=S_{n}$. For $(\pi_1,\cdots,\pi_m) \in Com (S_n,m)$, let
$\#(\pi_1,\cdots,\pi_m)$ be the number of connected components in the graph on
vertex set $\{1,\cdots,n\}$ defined by connecting the vertices according to the
permutations $\pi_1,\cdots,\pi_m$. For instance, $\#(\pi_1)$ is the number of
cycles of $\pi$. The main result of this note is the following theorem.

\begin{theorem} \label{mainresult} Let $\chi $ denote the (ordinary) Euler
characteristic of $M$. The generating function for the orbifold Euler
characteristic $\chi _{m}(M^{n},S_{n}) $ satisfies the following formulas:
\begin{eqnarray} \label{eqn: 1st main equality}
\sum_{n=0}^{\infty }\chi _{m}(M^{n},S_{n})u^{n}&=&\sum_{n=0}^{\infty
}\frac{u^{n}}{n!}\sum_{\pi _{1},\ldots,\pi _{m}\in Com (S_{n},m)}\chi ^{\#
(\pi _{1},\ldots,\pi _{m})}\\ \label{eqn: 2nd main equality}
&=&\left(\sum_{n=0}^{\infty }|Com (S_{n},m)|\frac{u^{n}}{n!} \right)^{\chi
}\\ \label{eqn: 3rd main equality}
&=&\prod_{i_1,\cdots,i_{m-1} = 1}^{\infty} ({1-u^{i_1 \cdots
i_{m-1}}})^{-\chi i_1^{m-2}i_2^{m-3} \cdots i_{m-2}}.
\end{eqnarray}
\end{theorem}

{\bf Remarks:} We will show that Equation \ref{eqn: 1st main equality}
follows directly from the definitions and a straightforward geometric
argument. Equation \ref{eqn: 2nd main equality} is proved in
Lemma \ref{t=1} and shows that it suffices to prove Equation
\ref{eqn: 3rd main equality} in the case $\chi =1$. Our main result then
should be regarded as Equation \ref{eqn: 3rd main equality} which in light
of Equation \ref{eqn: 2nd main equality} gives a generating function for
the number
of commuting $m$-tuples in $S_{n}$. Note also that for $m=1$ Equation
\ref{eqn: 3rd main equality} is Macdonald's formula $(1-u)^{-\chi}$ for the
Euler
characteristic of a symmetric product and for $m=2$ Equation \ref{eqn: 3rd
main equality} is Hirzebruch and H\"ofer's formula for the string-theoretic
orbifold Euler characteristic of a symmetric product.

To prove Equation \ref{eqn: 1st main equality} it suffices to see that
$$\chi (M^{(\pi _{1},\ldots,\pi _{m})})=(\chi (M))^{\# (\pi _{1},\ldots,\pi
_{m})}.$$
Partition $\{1,\ldots,n\}$ into disjoint subsets $I_{1},\ldots,I_{\#
(\pi _{1},\ldots,\pi _{m})}$ according to the
components of the graph associated to $(\pi _{1},\ldots,\pi _{m})$. Then
the small diagonal in the product $\prod_{i\in I_{j}}M_{i}$ is fixed by
$(\pi _{1},\ldots,\pi _{m})$ and is homeomorphic to $M$. The full fixed set
of $(\pi _{1},\ldots,\pi _{m}) $ is then the product of all the small
diagonals in the subproducts associated to the $I_{j}$'s. By the
multiplicative properties of Euler characteristic we see that $\chi
(M^{(\pi _{1},\ldots,\pi _{m})})=(\chi (M))^{\# (\pi _{1},\ldots,\pi
_{m})}.$

\begin{lemma} \label{t=1} For $\chi $ a natural number,

\[ \sum_{n=0}^{\infty} \frac{u^n}{n!} \sum_{\pi_1,\cdots,\pi_m \in Com
(S_{n},m)} \chi ^{\#(\pi_1,\cdots,\pi_m)} =
\left( \sum_{n=0}^{\infty} \frac{u^n |Com(S_n,m)|}{n!} \right) ^\chi  \]

\end{lemma}

\begin{proof}
        It suffices to show that an ordered $m$-tuple $(\pi_1,\cdots,\pi_m)$ of
mutually
commuting elements of $S_n$ contributes equally to the coefficient of
$\frac{u^n}{n!}$ on both sides of the equation. The contribution to this
coefficient on the left-hand side is $\chi ^{\#(\pi_1,\cdots,\pi_m)}$.

        The right hand side can be rewritten as

\[ \sum_{n=0}^{\infty} \frac{u^n}{n!} \sum_{n_1,\cdots,n_\chi : \sum n_i=n} {n
\choose n_1,\cdots,n_\chi} |Com(S_{n_1},m)| \cdots |Com(S_{n_\chi},m)|. \]

        Observe that ${n \choose n_1,\cdots,n_\chi} |Com(S_{n_1},m)| \cdots
|Com(S_{n_\chi},m)|$
is the number of ways of decomposing the vertex set $\{1,\cdots,n\}$ into
$\chi$
ordered subsets $S_1,\cdots,S_\chi$ of sizes $n_1,\cdots,n_\chi$ and
defining an
ordered
$m$-tuple of mutually commuting elements of $S_{n_i}$ on each subset. Gluing
these together defines an ordered $m$-tuple of mutually commuting elements of
$S_n$. Note that the $m$-tuple $(\pi_1,\cdots,\pi_m)$ arises in
$\chi^{\#(\pi_1,\cdots,\pi_m)}$ ways, because each of the
$\#(\pi_1,\cdots,\pi_m)$
connected components of the graph corresponding to $(\pi_1,\cdots,\pi_m)$ could
have come from any of the $\chi$ subsets $S_1,\cdots,S_\chi$.
\end{proof}

        Let us now recall some facts about wreath products of groups. All of
this can
be found in Sections 4.1 and 4.2 of James and Kerber \cite{James}. Given a
group
$G$, the wreath product $G Wr S_n$ is defined as a set by
$(g_1,\cdots,g_n;\pi)$
where $g_i \in G$ and $\pi \in S_n$. Letting permutations act on the right, the
group multiplication is defined by:

\[ (g_1,\cdots,g_n;\pi) (h_1,\cdots,h_n;\tau) = (g_1 h_{(1)\pi^{-1}}, \cdots,
g_n
h_{(n)\pi^{-1}};\pi \tau) \]

        Furthermore, the conjugacy classes of $G Wr S_n$ are parameterized as
follows. Let $Cl_1,\cdots,Cl_i$ be the conjugacy classes of $G$. Then the
conjugacy classes of $G Wr S_n$ correspond to arrays $(M_{j,k})$ satisfying the
properties:

\begin{enumerate}

\item $M_{j,k}=0$ if $j>i$
\item $\sum_{j,k} k M_{j,k}=n$

\end{enumerate}

        The correspondence can be made explicit. For $(g_1,\cdots,g_n;\pi) \in
G
Wr S_n$, let $M_{j,k}$ be the number of $k$-cycles of $\pi$ such that
multiplying
the $k$
$g_i$ whose subscripts lies in the $k$-cycle gives an element of $G$ belonging
to
the conjugacy class $Cl_j$ of $G$. The matrix so-defined clearly satisfies the
above two conditions.

        Lemma \ref{structure} is a key ingredient of this paper. It says that
centralizers of elements of wreath products can be expressed in terms
of
wreath products; this will lead to an inductive proof of Theorem
\ref{mainresult}.

\begin{lemma} \label{structure} Let $C_i$ denote a cyclic group of order $i$.
Then the centralizer in $C_i Wr S_n$ of an element in the conjugacy class
corresponding to the data $M_{j,k}$ is isomorphic to the direct product

\[ \prod_{j,k} C_{ik} Wr S_{M_{j,k}} \]

\end{lemma}

\begin{proof}
        To start, let us construct an element $(g_1,\cdots,g_n;\pi)$ of $C_i Wr
S_n$ with
conjugacy class data $M_{j,k}$. This can be done as follows:

\begin{enumerate}

\item Pick $\pi$ to be any permutation with $\sum_j M_{j,k}$ $k$-cycles

\item For each $j$ choose $M_{j,k}$ of the $k$-cycles of $\pi$ and think of
them
as $k$-cycles of type $j$

\item Assign (in any order) to the $g_i$ whose subscripts are contained in a
$k$-cycle of type $j$ of $\pi$ the values $(c_j,1,\cdots,1)$ where $c_j$ is an
element in the $jth$ conjugacy class of the group $C_i$

\end{enumerate}

        To describe the centralizer of this element $(g_1,\cdots,g_n;\pi)$,
note that
conjugation in $G Wr S_n$ works as

\begin{eqnarray*}
& & (h_1,\cdots,h_n;\tau) (g_1,\cdots,g_n;\pi) (h_{(1)\tau}^{-1},
\cdots,h_{(n)\tau}^{-1};\tau^{-1})\\
& = &(h_1 g_{(1)\tau^{-1}} h_{(1)\tau \pi^{-1} \tau^{-1}}^{-1},\cdots; \tau \pi
\tau^{-1})
\end{eqnarray*}

        It is easy to see that if $(h_1,\cdots,h_n;\tau)$ commutes with
$(g_1,\cdots,g_n;\pi)$, then $\tau$ operates on the $M_{j,k}$ $k$-cycles of
$\pi$
of type $j$ by first permuting these cycles amongst themselves and then
performing some power of a cyclic shift within each cycle. Further, among the
$h_i$ whose subscripts lie in a $k$-cycle of $\pi$ of type $j$ exactly one can
be
chosen arbitrarily in $C_i$--the other $h$'s with subscripts in that $k$-cycle
then have determined values.

        The direct product assertion of the theorem is then easily checked; the
only
non-trivial part is to see the copy of $C_{ik} Wr S_{M_{j,k}}$. Here the
$S_{M_{j,k}}$ permutes the $M_{j,k}$ $k$-cycles of type $j$, and the generator
of
the $C_{ik}$ corresponds to having $\tau$ cyclically permuting within the $k$
cycle and having the $h$'s with subscripts in that $k$-cycle equal to
$\{c_j,1,\cdots,1\}$, where $c_j$ is a generator of $C_i$.
\end{proof}

        With these preliminaries in hand, induction can be used to prove the
following
result. Note that by Lemma \ref{t=1}, only the $i=1$ case of Theorem
\ref{induction} is needed to prove the main result of this paper, Theorem
\ref{mainresult}. However, the stronger statement (general $i$) in Theorem
\ref{induction} makes the induction work by making the induction hypothesis
stronger.

\begin{theorem} \label{induction} For $m \geq 2$,

\[ \sum_{n=0}^{\infty} \frac{u^n |Com(C_i Wr S_n,m)|}{|C_i Wr S_n|} =
\prod_{i_1,\cdots,i_{m-1}=1}^{\infty} (\frac{1}{1-u^{i_1 \cdots
i_{m-1}}})^{i^{m-1} i_1^{m-2} i_2^{m-3} \cdots i_{m-2}} \]

\end{theorem}

\begin{proof}
        The proof proceeds by induction on $m$. We use the notation that if
$\lambda$ denotes a conjugacy class of a group
$G$, then $C_G(\lambda)$ is the centralizer in $G$ of some element of $\lambda$
(hence $C_G(\lambda)$ is well defined up to isomorphism). For the base case
$m=2$
observe that

\begin{eqnarray*}
& & \sum_{n=0}^{\infty} \frac{u^n |Com(C_i Wr S_n,2)|}{|C_i Wr S_n|}\\
& = & \sum_{n=0}^{\infty} \frac{u^n}{|C_i Wr S_n|} \sum_{(M_{j,k}): 1 \leq j
\leq
i \atop \ \ \ \sum_{j,k} k M_{j,k}=n} \frac{|C_i Wr S_n|}{|C_{C_i Wr
S_n}(M_{j,k})|} |C_{C_i Wr S_n}(M_{j,k})|\\
& = &  \sum_{n=0}^{\infty} u^n \sum_{(M_{j,k}): 1 \leq j \leq
i \atop \ \ \ \sum_{j,k} k M_{j,k}=n} 1 \\
& = & \prod_{i_1=1}^{\infty} (\frac{1}{1-u^{i_1}})^i
\end{eqnarray*}

        For the induction step, the parameterization of conjugacy classes of
wreath products and Lemma \ref{structure} imply that

\begin{eqnarray*}
& & \sum_{n=0}^{\infty} \frac{u^n |Com(C_i Wr S_n,m)|}{|C_i Wr S_n|}\\
& = & \sum_{n=0}^{\infty} \frac{u^n}{|C_i Wr S_n|} \sum_{(M_{j,k}): 1 \leq j
\leq
i \atop \ \ \ \sum_{j,k} k M_{j,k}=n} \frac{|C_i Wr S_n|}{|C_{C_i Wr
S_n}(M_{j,k})|} |Com(C_{C_i Wr S_n}(M_{j,k}),m-1)|\\
& = & [\prod_{k=1}^{\infty} \sum_{a=0}^{\infty} \frac{u^{ka}|Com(C_{ik} Wr
S_a,m-1)|}{|C_{ik} Wr S_a|}]^i\\
& = & [\prod_{k=1}^{\infty} \prod_{i_2,\cdots,i_{m-1}=1}^{\infty}
(\frac{1}{1-u^{ki_2
\cdots i_{m-1}}})^{(ik)^{m-2} i_2^{m-3} \cdots i_{m-2}}]^i\\
& = & \prod_{i_1,\cdots,i_{m-1}=1}^{\infty} (\frac{1}{1-u^{i_1 \cdots
i_{m-1}}})^{i^{m-1} i_1^{m-2} i_2^{m-3} \cdots i_{m-2}}
\end{eqnarray*}

\end{proof}

\textsc{\\
Department of Mathematics\\
University of California\\
Berkeley, CA 94720}

\textsc{\\
Department of Mathematics\\
Dartmouth College\\
Hanover, NH 03755}

\texttt{\\
jbryan{@}math.berkeley.edu\\
jason.e.fulman{@}dartmouth.edu}

\end{document}